\documentclass[11pt]{amsart}

\usepackage[square,compress,comma, numbers,sort]{natbib}
\usepackage[colorlinks=true, citecolor=red, linkcolor=blue]{hyperref}
\usepackage{amsfonts,mathtools}

\allowdisplaybreaks[4]

\usepackage{amsmath}
\usepackage{bbm}
\usepackage{amssymb}
\usepackage{mathtools}

\usepackage{color}

\definecolor{c20}{rgb}{0.,0.7,0.}
\definecolor{c30}{rgb}{0.,0.,1.}
\definecolor{c40}{rgb}{1,0.1,0.7}
\definecolor{c50}{rgb}{1,0,0}
\definecolor{c60}{rgb}{1,0.9,0.1}

\newcommand{\abs}[1]{\left\lvert #1 \right\rvert}

\newcommand{\pk}[1]{\mathbb{P} \left\{ #1 \right \} }

\newcommand{\R}{\mathbb{R}}

\newcommand{\inr}{\in \R}

\newcommand{\ldot}{,\ldots,}

\newcommand{\BQN}{\begin{eqnarray}}
\newcommand{\EQN}{\end{eqnarray}}
\newcommand{\BQNY}{\begin{eqnarray*}}
\newcommand{\EQNY}{\end{eqnarray*}}

\newcommand{\BS}{\begin{sat}}
\newcommand{\ES}{\end{sat}}
\newcommand{\BT}{\begin{theo}}
\newcommand{\ET}{\end{theo}}
\newcommand{\BK}{\begin{korr}}
\newcommand{\EK}{\end{korr}}

\newcommand{\BD}{\begin{de}}
\newcommand{\ED}{\end{de}}

\newcommand{\BIT}{\begin{itemize}}
\newcommand{\EIT}{\end{itemize}}
\newcommand{\BDI}{\begin{description}}
\newcommand{\EDI}{\end{description}}

\newcommand{\BRM}{\begin{remarks}}
\newcommand{\ERM}{\end{remarks}}

\newcommand{\BEL}{\begin{lem}}
\newcommand{\EEL}{\end{lem}}

\newtheorem{theo}{Theorem}[section]
\newtheorem{sat}[theo]{Proposition}
\newtheorem{de}[theo]{Definition}
\newtheorem{lem}[theo]{Lemma}

\newtheorem{korr}[theo]{Corollary}
\newtheorem{remark}[theo]{Remark}
\newtheorem{remarks}[theo]{Remarks}

\newcommand{\nelem}[1]{{Lemma \ref{#1}}}
\newcommand{\nedef}[1]{{Definition \ref{#1}}}

\newcommand{\netheo}[1]{{Theorem \ref{#1}}}

\newcommand{\prooftheo}[1]{ \textsc{\bf Proof of Theorem} \ref{#1}:}

\newcommand{\prooflem}[1]{\textsc{\bf Proof of Lemma} \ref{#1}:}

\newcommand{\COM}[1]{}

\def\td{\text{\rm d}}

\newcommand{\QED}{\hfill $\Box$}

\newcommand{\kb}[1]{\boldsymbol{#1}}
\newcommand{\vk}[1]{\kb{#1}}

\def\bqny#1{{\begin{eqnarray*} #1 \end{eqnarray*}}}

\begin{document}

\title{Uniform   bounds for ruin probability in Multidimensional Risk Model}

\author{Nikolai Kriukov }
\address{Nikolai Kriukov, Department of Actuarial Science, %\\Faculty of Business and Economics\\
University of Lausanne,\\
UNIL-Dorigny, 1015 Lausanne, Switzerland
}
\email{Nikolai.Kriukov@unil.ch}

\bigskip

\date{\today}
\maketitle

{\bf Abstract:} In this paper we consider some generalizations of the classical $d$-dimensional Brownian risk model. This contribution derives some non-asymptotic bounds for simultaneous ruin probabilities of interest. In addition, we obtain non-asymptotic bounds also for the case of general trend functions and convolutions of our original risk model.


  
\section{Introduction and first Result} 
Let $\vk B(t),t\ge 0$ be a $d$-dimensional Brownian motion with  independent standard Brownian motion components and set $\vk Z(t)= A \vk B(t), t\ge 0$ 
with $A$ a $d\times d$ real non-singular matrix. The recent contribution \cite{KWW} derived the following remarkable inequality
\BQN
1 \le \frac{\pk{ \exists t\in [0,T]:\vk{Z}(t) \ge \vk{b}}}{\pk{\vk{Z}(T)\ge \vk{b}} } &\le& K(T), \quad   K(T)= \frac{1}
{\pk{\vk{Z}(T)\ge \vk 0} }
\label{Korsh0}
\EQN
valid for all $\vk b \inr^d, T>0$. In our notation bold symbols are column vectors with $d$ rows and all operations are meant component-wise, for instance $\vk x \ge \vk 0$ means $x_i \ge 0$ for all $i\le d$ with $\vk 0= (0 \ldot 0)\inr^d$.  \\
The special and crucial feature of  \eqref{Korsh0} is that the bounds are uniform with respect to $\vk b$. Moreover,  if at least one component of $\vk b$ tends to infinity, then $\pk{ \exists t\in [0,T]:\vk{Z}(t) \ge \vk{b}}$ 
can be accurately approximated  (up to some constant) by the survival probability $\pk{\vk{Z}(T)\ge \vk{b}}$.\\
 Inequality \eqref{Korsh0} has been crucial in the context of Shepp-statistics investigated in \cite{KWW}. It is also of great importance in the investigation of simultaneous ruin probabilities in vector-valued risk models (see \cite{ZKE,parisian,doi:10.1080/03461238.2021.1902853}). Specifically, consider 
the multidimensional risk model
$$ \vk R(t,u)=\vk a u-\vk X(t), \  \vk X(t)= \vk Z(t)- \vk c t$$
for some vectors $\vk a,\vk c\in\R^d$ and $\vk Z(t),t\ge 0$ defined above. Typically,   $\vk R$ models the surplus of all $d$-portfolios of an insurance company, where $a_i u, u>0$ plays the role of the initial capital. Here the component $Z_i$ models the accumulated claim amount up to time $t$ and $c_i t$ is the premium income for the $i$th portfolio.    
 
Given a positive integer $k\le d$, of interest is the calculation of the $k$-th simultaneous ruin probability, i.e.,  at least $k$ out of $d$ portfolios are ruined on a given time interval $[0,T]$ with $T$ possibly also infinite. That ruin probability can be written as  
$$\pk{\exists_{t\in[0,T]}:\vk Z(t)-\vk c t\in u \vk S}, \qquad u>0,$$
where 
$$\vk S:=\bigcup_{\substack{I\subset\{1\ldots d\}\\\abs{I}=k}}\vk S_{I},\qquad \vk S_{I}=\{\vk x\in\R^d:~\forall i\in I,~x_i>a_i\}.$$
The particular case $\vk Z(t)= A \vk B(t), t\ge 0$ with $A$ a $d\times d$ non-singular matrix is of special importance for insurance risk models, see 
e.g., \cite{delsing2018asymptotics}. Clearly, this instance is  also of great importance in statistics and probability given the central role of the $\R^d$-valued Brownian motion as a natural limiting process\\
 In \cite{pandemic} it has been shown that 
\eqref{Korsh0} can be extended for this risk model, i.e., for all $u,T>0$
\BQN
1 \le \frac{\pk{\exists_{t\in [0,T]}: \vk X(t) \in u \vk S} }{  \pk{ \vk X(T)\in u \vk S} } &\le& K_{ \vk S}(T), \  \vk X(t)=\vk Z(t)-\vk c( t), 
\label{Korsh1}
\EQN
with $\vk c(t)= \vk c t, t\ge 0$ and some known constant $K_{\vk S}(T)>0$. 
Again the bounds are uniform with respect to $u$.\\
 It is clear that the inequality \eqref{Korsh1} does not hold for an arbitrary set $\vk S\subset\R^d$. Since Brownian motion has almost surely continuous sample paths, , if it hits some closed set $\vk S$, it definitely hits its boundary. Hence in the following special case, for all $u$ positive we have
\bqny{
\{\exists t\in[0,T]:Z(t)\geq u\}=\{\exists t\in[0,T]:Z(t)=u\}.
} 
Hence, taking $\vk S=\{\vk x\in\R^d: x_1=1\}$ and $\vk c(t)=\vk 0$ we have that
\bqny{
\pk{\exists_{t\in [0,T]}: \vk X(t) \in u \vk S}&=&\pk{\exists_{t\in [0,T]}:  X_1(t) > u}\geq\pk{X_1(T) > u}>0,  \\
\pk{ \vk X(T)\in u \vk S}&=&0,
}
and \eqref{Korsh1} does not hold.
\\
Therefore hereafter we shall consider only closed sets S described as follows:

\begin{de}\label{cone_def} Let $\vk X$ and $\vk Z$ are as defined above. The closed Borel set $\vk{S}\subset \R^d$ satisfies the cone condition with respect to the the vector-valued process $\vk X$ if there exists a strictly positive function $\varepsilon_{\vk{S}}(t), t> 0$ such that for any point $\vk x\in \vk{S}$  and any $t>0$ there exists a Borel set  $\vk V_{\vk x}\subset \vk{S}$ that contains  $\vk{x}$ and not depending on t, satisfying $\vk V_{\vk x}-\vk x \subset C \left(\vk V_{\vk x}-\vk x\right)$ for all $C>1$ and  $\pk{ \vk Z(t) \in \vk V_x-\vk x} \geq\varepsilon_{\vk S}(t)$.
	\end{de}

It is of interest to consider a general trend function in \eqref{Korsh1}. We consider below a large class of trend functions which is tractable  if  $\vk Z$ has self-similar coordinates with index $\alpha>0$. This is in particular the case  when $\vk Z= A \vk B $.
\begin{de}
A continuous measurable vector-valued function $\vk c:[0,+\infty)\to\R^d$ belongs to $RV_{t_0}(\alpha)$ for some $\alpha>0$, $t_0\in[0,T]$ if for some $M>0$, all $i\in\{1\ldots d\}$, $t\in[0,T]$ 
$$|c_i(t)-c_i(t_0)|\le M|t-t_0|^{\alpha}.$$
\end{de}

 We state next our first result. Below  $\vk F:\R^d\to\R^d$ growing means that  for any $\vk x,\vk y\in\R^d$ such that for all $i\in\{1\ldot d\}$ $x_i\geq y_i$ we have that $F_i(\vk x)\geq F_i(\vk y)$ holds for all $i\in\{1\ldot d\}$.

\begin{theo}\label{ConeSet}
If  $\vk S \subset \R^d$ satisfies the cone condition with respect to the process $ \vk Z=A \vk B$ such that  $\vk 0\not\in\vk S$  and $\vk c \in  RV_T(1/2),$ then for all constants $T>0, u>1$ the inequality 
\eqref{Korsh1} holds %with $\vk Z=  A  \vk B$ 
with 
	\bqny{
K(T)= %		\pk{\exists_{t\in[0,T]}:\vk W(t)-\vk c(t)\in u\vk S}\leq
\frac{2^{d/2}}{\mathfrak{C}(T)\varepsilon_{\vk S}(T)}, \quad 
		\mathfrak{C}(T)=\inf_{t\in [0,T)}e^{-T\left(\frac{\vk c(T)-\vk c(t)}{\sqrt{T-t}}\right)^\top\Sigma^{-1}\left(\frac{\vk c(T)-\vk c(t)}{\sqrt{T-t}}\right)}>0,
	}
where $\Sigma$ is the covariance matrix of $\vk Z(T)$.
 In particular, for any growing function $\vk F:\R^d\to\R^d$
\bqny{
	\pk{\exists_{t\in [0,T] }: \vk F(\vk Z(t)-\vk c(t))>u\vk a}\leq C_T\pk{\vk F(\vk Z(T)-\vk c(T))>u\vk a}
}
we have $\vk a \in \R^d\setminus (-\infty,0]^d$,  $u>1$ and  some constant $C_T$ which does not depend on $u$.
\end{theo}

If  $\vk Z$ is a given separable random field, it is of interest to determine conditions such that \eqref{Korsh1} can be extended to  
\BQN
1 \le \frac{\pk{\exists_{ \vk t \in \mathbb{T} }:\vk Z(\vk t)-\vk c(\vk t)\in u \vk S} }{  \pk{ \vk Z(\vk T)-\vk c(\vk T)\in u \vk S} } &\le& K_{ \vk S}(\vk T), 
\label{Korsh2}
\EQN 
where $\mathbb{T} = [0,T_1]\times\ldots\times[0,T_n]$  and $\vk T=(T_1 \ldot T_n)$ has positive components. For the case $\vk Z(\vk t)= \sum_{i=1}^n \vk Z_i(t_i)$ where $\vk Z_i$'s are independent copies of $\vk Z$,  and $\vk c(\vk t)=0$  the result \eqref{Korsh2} was shown in  \cite{KWW}[Thm 1.1] for some special set $\vk S$. For more general set $\vk S$ we have the following result:
 
\begin{theo}\label{convshift}
If  $\vk S\subset \R^d$ satisfies the cone condition with respect to $\vk Z$, $\vk 0\not\in\vk S$  and  all $\vk c_i \in  RV_T(1/2),$ then for all constants $T_1\ldot T_n>0, u>1$ the inequality 
\eqref{Korsh2} holds with $\vk Z(\vk t)= \sum_{k=1}^{n}\vk Z_k(t_k)$ and $\vk c(\vk t)=\sum_{k=1}^{n}\vk c_k(t_k)$ with 
	\bqny{
K_{\vk S}(\vk T)=\prod_{k=1}^{n}\frac{2^{d/2}}{\mathfrak{C}_k(T_k)\varepsilon_{\vk S}(T_k)}, \quad 
		\mathfrak{C}_k(T_k)=\inf_{t\in[0,T_k)}e^{-T_k\left(\frac{\vk c_k(T_k)-\vk c_k(t)}{\sqrt{T_k-t}}\right)^\top\Sigma^{-1}(T_k)\left(\frac{\vk c_k(T_k)-\vk c_k(t)}{\sqrt{T_k-t}}\right)}>0,
	}
where $\varepsilon_{\vk S}$ is any function satisfies the claims of Definition 1.1.
\end{theo}

\section{Discussion}
In this section as in Introduction
we consider  first
$$ \vk Z(t)=A\vk B(t),\qquad t\geq 0$$
with $A$ non-singular and $\vk B$ a $d$-dimensional Brownian motion with independent components. 
 We shall discuss next the generalisation of the upper bound \eqref{Korsh1} for various special cases.

\subsection{Order statistics}
The classical multidimensional Brownian motion risk model (see \cite{delsing2018asymptotics}) is formulated in terms of some risk process $\vk R$ specified by 
$$ \vk R(t,u)=\vk a u-\vk Z(t)+\vk c t$$
for some vectors $\vk a,\vk c\in\R^d$. 
We are interested in the finite-time simultaneous ruin probability for $k$ out of $d$ portfolios,  i.e. the probability that at least $k$ portfolios are ruined. In other words, we are investigating the probability 
$$\pk{\exists_{t\in[0,T]},\exists_{\mathcal{I}\subset\{1,\ldots,d\}}:\abs{\mathcal{I}}=k,~\forall i\in\mathcal{I}~Z_i(t)-c_it\geq a_i u},$$
which in view of our previous notation reads
$$\pk{\exists_{t\in[0,T]}:\vk Z(t)-\vk c t\in \vk S_u}, \qquad u>0,$$
where
$$\vk S_u:=\bigcup_{\substack{I\subset\{1\ldots d\}\\\abs{I}=k}}\vk S_{I,u},\qquad \vk S_{I,u}=\{\vk x\in\R^d:~\forall i\in I~x_i\geq a_iu\}.$$
Asymptotic approximations of such probability was already obtained in \cite{pandemic}. Now we want to derive a uniform non-asymptotic bound based on our previous findings. It is clear that all sets $\vk S_{I,u}$ satisfy the cone condition with respect to the process $\vk Z$. Thus, $\vk S_u$ also satisfies the cone condition with respect to the process $\vk Z$, hence we can use \netheo{ConeSet} and write for some positive constant $C$
$$ \pk{\vk Z(T)-\vk cT\in\vk S_u}\leq \pk{\exists_{t\in[0,T]}:\vk Z(t)-\vk c t\in\vk S_u}\leq C\pk{\vk Z(T)-\vk cT\in\vk S_u}.$$

\subsection{Fractional Brownian motion}
Consider next the 1-dimensional risk model
$$R(u,t)=u-B_H(t)+ct, \, t >0,$$ 
where $B_H(t)$ is a standard fractional Brownian motion with zero mean and variance function  $\abs{t}^{2H}$, $H\in (0,1]$. 
We are interested in the calculation of the finite-time ruin probability for  given $T>0$. The inequalities below have already been shown in \cite{ZKE}. We retrieve them using  our findings. Namely, by the Slepian inequality, we can write for $H>\frac{1}{2}$ and W a standard Brownian motion
\bqny{
\pk{\exists_{t\in[0,T]}R(u,t)\leq 0}&\leq&\pk{\exists_{t\in[0,T]}W\left(t^{2H}\right)-ct\geq u}\\
&=&\pk{\exists_{t\in[0,T^{2H}]}W\left(t\right)-ct^{1/2H}\geq u}\\ 
&=&\pk{\exists_{t\in[0,1]}W\left(T^{2H}t\right)-cTt^{1/2H}\geq u}\\
&=&\pk{\exists_{t\in[0,1]}W\left(t\right)-cT^{1-H}t^{1/2H}\geq u/T^H}.
}
Since $cT^{1-H}t^{1/2H}\in RV_{1}(1/2)$, using \netheo{ConeSet}, for some positive constant $C$ we can write
\bqny{
\pk{\exists_{t\in[0,1]}W\left(t\right)-cT^{1-H}t^{1/2H}\geq u/T^H}&\leq&C\pk{W\left(1\right)-cT^{1-H}\geq u/T^H}\\
&=&C\pk{W\left(T^{2H}\right)-cT\geq u}\\
&=&C\pk{R(u,T)\leq 0}.
}

%\subsection{Convolution of fractional Brownian motions}
The above can be extended considering the convolution of $n$ independent one-dimensional fractional Brownian motions $B_{H_i}(t), t>0, i\le n$. Let $H_i>1/2$ and define the risk processes  $$R_i(u,t)=u/n-B_{H_i}(t)+c_it,\,  i\le n. $$ 
Consider the convolution of processes $R_i(u,t)$. Using Slepian inequality, for all $H_i>\frac{1}{2}$ we can write 
\bqny{
\pk{\exists_{\vk t\in\prod\limits_{i=1}^{n}[0,T_i]}\sum_{i=1}^{n}R_i(u,t_i)\leq 0}&\leq&\pk{\exists_{\vk t\in\prod\limits_{i=1}^{n}[0,T_i]}\sum_{i=1}^{n}W_{i}\left(t^{2H_i}\right)-c_it_i\geq u}\\
&=&\pk{\exists_{\vk t\in\prod\limits_{i=1}^{n}[0,T_i^{2H_i}]}\sum_{i=1}^nW_{i}\left(t\right)-c_it_i^{1/2H_i}\geq u}.
}
Here $W_i$ stands for an independent copy of Brownian motion.
As $c_it^{1/2H_i}\in RV_{T_i}(1/2,1)$, using \netheo{convshift}, for some positive constant $C$ we can write
\bqny{
\pk{\exists_{\vk t\in\prod\limits_{i=1}^{n}[0,T_i^{2H_i}]}\sum_{i=1}^nW_{i}\left(t\right)-c_it_i^{1/2H_i}\geq u}&\leq&C\pk{\sum_{i=1}^{n}W_{i}\left(T_i^{2H_i}\right)-c_iT_i\geq u}\\
&=&C\pk{\sum_{i=1}^{n}B_{H_i}\left(T_i\right)-c_iT_i\geq u}\\
&=&C\pk{\sum_{i=1}^{n}R_i(u,T_i)\leq 0}.
}

\section{Vector-valued time-transform}
Finally, we discuss some extensions of \eqref{Korsh1} under different time transformations. We use the notation from Section 2 and define the following time transform. Let $\vk f(t):[0,+\infty)\in\R^d$ be a growing vector-valued function and define
\bqny{
\vk Z(\vk f(t))=(Z_1(f_1(t)),\ldots, Z_d(f_d(t)))^\top.
} 
Hence $\vk f(t)$ can be considered as a generalised transformation of time. 
\begin{theo}\label{TimeTransform}
Let $\vk c(t),\vk f(t):[0,T]\to\R^d$ be given. Suppose that all $f_i(t)$'s are continuous, strictly growing and for all $i\in\{1\ldot d\}$ we have $f_i(0)=0$ and function $\delta_i(t)=\frac{f_i(T)-f_i(t)}{f_1(T)-f_1(t)}$ has a positive finite limit as $t\to T$. Let also $\abs{c_i(T)-c_i(t)}< M\sqrt{f_1(T)-f_1(t)}$ for all $t\in[0,T]$, all $i\in\{1\ldot d\}$, some $M>0$, and $\vk S$ satisfies the cone condition with respect to the process $\vk Z$. If $\vk 0\not\in\vk S$, then for all constants $T>0, u>1$ the inequality 
\eqref{Korsh1} holds with $\vk X(t)= \vk Z(\vk f(t))$ and 
	\bqny{
K^*(T)= \frac{(2f_1(T))^{d/2}}{\mathfrak{C}(T)\bar{\varepsilon}_{\vk S}}, \quad 
		\mathfrak{C}(T)=\inf_{t\in [0,T)}e^{-\left(\frac{\vk c(T_k)-\vk c(t)}{\sqrt{f_1(T)-f_1(t)}}\right)^\top\Sigma^{-1}(\vk \delta(t))\left(\frac{\vk c(T_k)-\vk c(t)}{\sqrt{f_1(T)-f_1(t)}}\right)}>0,
	}
where
\bqny{
\bar\varepsilon_{\vk S}=\left(\frac{\inf\limits_{\substack{i\in\{1\ldot d\}\\t\in [0,T]}}\delta_i(t)}{\sup\limits_{\substack{i\in\{1\ldot d\}\\t\in [0,T]}}\delta_i(t)}\right)^{d/2}\varepsilon_{\vk S}\left(\inf\limits_{\substack{i\in\{1\ldot d\}\\t\in [0,T]}}\delta_i(t)\right)>0.
}
\end{theo}

\begin{remark}
The function $\vk f$ in \netheo{TimeTransform} may also be an almost surely growing stochastic process, independent of $\vk Z$, satisfying 
\bqny{
& &\max_{i\in\{1,\ldots,d\}}f_i(T)<F,\qquad \max_{i\in\{1,\ldots,d\}}\sup_{t\in [0,T)}\abs{\frac{c_i(T_k)-c_i(t)}{\sqrt{f_1(T)-f_1(t)}}}<M,\\
& &\delta<\inf\limits_{\substack{i\in\{1\ldot d\}\\t\in [0,T]}}\delta_i(t)\leq \sup\limits_{\substack{i\in\{1\ldot d\}\\t\in [0,T]}}\delta_i(t)<\Delta,
} 
almost surely with some positive constants $F,M,\delta,\Delta$. In this case the inequality \eqref{Korsh1} holds with
\bqny{
K^*(T)= \frac{(2F)^{d/2}}{\mathfrak{C}(T)\bar{\varepsilon}_{\vk S}}, \quad 
		\mathfrak{C}(T)=\min_{\substack{\vk x\in [-M,M]^d\\ \vk t\in[\delta,\Delta]^d}}e^{-\vk x^\top\Sigma^{-1}(\vk t)\vk x}>0,
	}
and
\bqny{
\bar\varepsilon_{\vk S}=\left(\frac{\delta}{\Delta}\right)^{d/2}\varepsilon_{\vk S}\left(\delta\right)>0.
}
\end{remark}

We illustrate the above findings considering again 
%\subsection{Example. Independent fractional Brownian motions}
$d$ independent one-dimensional fractional Brownian motions $B_{H_i}(t), t>0$ with  Hurst parameters $H_i>\frac{1}{2}, i \le  d$. Define $d$ ruin portfolios $$R_i(u,t)=u-B_{H_i}(t)+c_it,$$ and we are interested in probability that all of them will be simultaneously ruined in  $[0,T]$.\\
 Using Gordon inequality (see\COM{e.g. \cite{MR800188} or} \cite[page~55]{MR1088478}), we obtain 
\bqny{
\pk{\exists_{t\in[0,T]}\forall_{i\in\{1 \ldot d\}} R_i(u,t)<0}&\leq&\pk{\exists_{t\in[0,T]} \forall_{i\in\{1 \ldot d\}}W_{i}\left(t^{2H_i}\right)-c_it>u}.
}
Where $B_i(t)$ are independent Brownian motions.
Since 
\bqny{
\lim_{t\to T}\frac{T^{2H_i}-t^{2H_i}}{T^{2H_1}-t^{2H_1}}=\frac{2H_i}{2H_1}\frac{T^{2H_i-1}}{T^{2H_1-1}}>0,
} 
using \netheo{TimeTransform}, for some positive constant $C$, which does not depend on $u$ we can write
\bqny{
\pk{\exists_{t\in[0,T]} \forall_{i\in\{1 \ldot d\}}W_{i}\left(t^{2H_i}\right)-c_it>u}&\leq&C\pk{\forall_{i\in\{1 \ldot d\}}W_{i}\left(T^{2H_i}\right)-c_iT>u}\\
&=&C\pk{\forall_{i\in\{1 \ldot d\}}B_{H_i}\left(T\right)-cT>u}\\
&=&C\pk{\forall_{i\in\{1 \ldot d\}}R_i(u,T)<0}.
}

\section{Proofs}
Let us note the following property of the function $\varepsilon_{\vk S}(t)$.
\begin{lem}\label{1} If set $\vk S$ satisfies the cone condition with respect to the process $\vk Z(t)$ with some function $\varepsilon_{\vk S}(t)$, then for any constant $u>1$ set $u\vk S$ also satisfies the cone condition with respect to the process $\vk Z(t)$, and for any function $\varepsilon_{\vk S}(t)$ exists a function $\varepsilon_{u\vk S}(t)$ such that $$\varepsilon_{u\vk S}(t)\geq \varepsilon_{\vk S}(t)$$
\end{lem}

\prooflem{1} Fix some $\vk x\in u\vk S$. Then we know that $\vk y=\vk x/u\in\vk S$. As $\vk S$ satisfies the cone condition with respect to the process $\vk Z(t)$, there exists some cone $V_{\vk y}\subset \vk S$ with vertex $\vk y$ such that $\pk{\vk Z(t)\in V_{\vk y}-\vk y}\geq \varepsilon_{\vk S}(t)$. Hence, $uV_{\vk y}\subset u\vk S$ for all $u>1$. Note that using the properties of cone
\bqny{
uV_{\vk y}=u(\vk y+(V_{\vk y}-\vk y))=\vk x+u(V_{\vk y}-\vk y)\supset \vk x+(V_{\vk y}-\vk y).
}
Hence, $\vk x+(V_{\vk y}-\vk y)\subset u\vk S$ is some cone with vertex $\vk x$, and 
\bqny{
\pk{\vk Z(t)\in uV_{\vk y}-\vk x}\ge\pk{\vk Z(t)\in V_{\vk y}-\vk y}\ge\varepsilon_{\vk S}(t).
}
\qed

\prooftheo{ConeSet} Consider the first inequality.
Define the following stopping moment
\bqny{
\tau=\inf\{t\in [0,T]:\vk Z(t)-\vk c(t)\in u\vk S \}.
}
According to the strong Markov property
\bqny{
& &\pk{\vk Z(T) -\vk c (T)\in u\vk S}\\
& &\qquad\quad=\int_{0}^{T}\int_{u \partial \vk S}\pk{\vk Z(\tau)-\vk c (\tau)\in\td\vk x,\tau\in \td t}\pk{\vk Z(T)-\vk c(T)\in u\vk S|\vk Z(t)-\vk c(t)=\vk x}.
}
Using \nelem{1}, $\vk uS$ satisfies the cone condition with respect to the process $\vk Z(t)$. Hence for all $\vk x\in u\vk S$, $t\in [0,T]$
\bqny{
\pk{\vk Z(T)-\vk c(T)\in u\vk S|\vk Z(t)-\vk c(t)=\vk x}&\geq&\pk{\vk Z(T)-\vk c(T)\in \vk V_{\vk x}|\vk Z(t)-\vk c(t)=\vk x}\\
&=&\pk{\vk Z(T-t)-(\vk c(T)-\vk c(t))\in \vk V_{\vk x}-\vk x}\\
&=&\pk{\vk Z(1)-(\vk c(T)-\vk c(t))/\sqrt{T-t}\in (\vk V_{\vk x}-\vk x)/\sqrt{T-t}}\\
&\geq&\pk{\sqrt{T}\vk Z(1)\in \vk V_{\vk x}-\vk x+\sqrt{T}(\vk c(T)-\vk c(t))/\sqrt{T-t}}\\
&=&\pk{\vk Z(T)\in \vk V_{\vk x}-\vk x+\sqrt{T}(\vk c(T)-\vk c(t))/\sqrt{T-t}}\\
&=&\int_{\vk V_{\vk x}-\vk x}\frac{1}{(2\pi)^\frac{d}{2}\sqrt{\abs{\Sigma}}}e^{-\frac{1}{2}(\tilde{\vk x}+\sqrt{T}\frac{\vk c(T)-\vk c(t)}{\sqrt{T-t}})^\top\Sigma^{-1}(\tilde{\vk x}+\sqrt{T}\frac{\vk c(T)-\vk c(t)}{\sqrt{T-t}})}\td\tilde{\vk x}\\
&\geq&\int_{\vk V_{\vk x}-\vk x}\frac{1}{(2\pi)^\frac{d}{2}\sqrt{\abs{\Sigma}}}e^{-T(\frac{\vk c(T)-\vk c(t)}{\sqrt{T-t}})^\top\Sigma^{-1}(\frac{\vk c(T)-\vk c(t)}{\sqrt{T-t}})}e^{-\frac{1}{2}(\sqrt{2}\tilde{\vk x})^\top\Sigma^{-1}(\sqrt{2}\tilde{\vk x})}\td\tilde{\vk x}\\
&=&e^{-T(\frac{\vk c(T)-\vk c(t)}{\sqrt{T-t}})^\top\Sigma^{-1}(\frac{\vk c(T)-\vk c(t)}{\sqrt{T-t}})}\frac{\pk{\vk Z(T)\in\sqrt{2}(\vk V_{\vk x}-\vk x)}}{2^{d/2}}\\
&\geq&\frac{1}{2^{d/2}}e^{-T(\frac{\vk c(T)-\vk c(t)}{\sqrt{T-t}})^\top\Sigma^{-1}(\frac{\vk c(T)-\vk c(t)}{\sqrt{T-t}})}\pk{\vk Z(T)\in\vk V_{\vk x}-\vk x}\\
&\geq&\frac{\mathfrak{C}(T)\varepsilon_{u \vk S}(T)}{2^{d/2}}\geq\frac{\mathfrak{C}(T)\varepsilon_{\vk S}(T)}{2^{d/2}},
}
where $\vk V_{\vk x}$ is the cone from \nedef{cone_def}. As the right part does not depend on $\vk x$ and $t$, we can write
\bqny{
\pk{\vk Z(T) -\vk c (T)\in u\vk S}&\geq&\frac{\mathfrak{C}(T)\varepsilon_{\vk S}(T)}{2^{d/2}}\int_{0}^{T}\int_{u\partial \vk S}\pk{\vk Z(\tau)-\vk c (\tau)\in\td\vk x,\tau\in \td t}\\
&=&\frac{\mathfrak{C}(T)\varepsilon_{\vk S}(T)}{2^{d/2}}\pk{\exists_{t\in [0,T]}\vk Z(t)-\vk c(t)\in u\vk S}.
}
Hence, the first inequality holds. Consider the second one. Define a set
\bqny{
\vk a^+=\{\vk x\in\R^d:\vk x\geq \vk a\}
}
and
\bqny{
\vk S_u=\frac{1}{u}\{\vk x\in\R^d: \vk F(x)\in u\vk a^+\}.
}
Set $\vk S_u$ satisfies the cone condition with respect to the process $\vk Z(t)$ for $V_{\vk x}=\vk x^+$, as for any $\vk y\geq \vk x\in\vk S_u$
\bqny{
\vk F(u\vk y)\geq\vk F(u\vk x)\geq u\vk a^+.
}
Consequently, $\vk y\in\vk S_u$, and 
\bqny{
\varepsilon_{\vk S_u}(t)=\pk{\vk X(t)\in\vk x^+-\vk x}=\pk{\vk X(t)\in[0,+\infty)^d}
}
does not depend on $u$. Applying the result above for the set $\vk S_u$ we obtain
\bqny{
\pk{\exists_{t\in [0,T]}:\vk X(t)\in u\vk S_u}\leq\frac{2^{d/2}\pk{\vk X(T)\in u\vk S}}{\mathfrak{C}(T)\varepsilon_{\vk S_u}(T)}=\frac{2^{d/2}\pk{\vk X(T)\in u\vk S_u}}{\mathfrak{C}(T)\pk{\vk X(T)\in[0,+\infty)^d}}.
}
As the event $\{\vk X(t)\in u\vk S_u\}$ is equal to the event $\{\vk F(\vk X(t)-\vk c(t))>u\vk a\}$, this completes the proof.
\QED

\prooftheo{convshift}
Define
\bqny{
\psi_k(\vk S):=\pk{\exists_{\vk t\in\mathbb{T}_k}:\sum_{i=1}^{k}(\vk Z_i(t_i)-\vk c_i(t_i))+\sum_{i=k+1}^{n}(\vk Z_i(T_i)-\vk c_i(T_i))\in\vk S},
}
where $\mathbb{T}_k=[0,T_1]\times\ldots\times[0,T_k]$. As in the previous section we are going to prove that the inequality
\bqny{
\psi_{k}(u\vk S)\leq\frac{2^{d/2}\psi_{k-1}(u \vk S)}{\varepsilon_{\vk S}(T_k)\mathfrak{C}_k(T_k)}
} 
takes place for any $k\in\{1\ldot n\}$. We can fix the trajectories of processes $\vk Z_i(t)$ called $\vk x_i(t)$, fix random vectors $\vk Z_i(T_i)$ called $\vk x_i$, and define the process
\bqny{
\vk Z^{*k}(t,\vk t^k)=\vk Z_k(t)-\vk c_k(t)+\sum_{i=1}^{k-1}(\vk x_i(t_i)-\vk c_i(t_i))+\sum_{i=k+1}^{n}(\vk x_i-\vk c_i(T_i)),
}
where $\vk t^k=(t_1\ldot t_{k-1})\in\mathbb{T}_{k-1}$. \\
Since $\vk Z_i$ are independent, it is enough to show that for every set of trajectories $\vk x_i(t)$ and points $\vk x_j$, the inequality
\bqny{
\psi^*(u\vk S)\leq\frac{2^{d/2}\nu(u\vk S)}{\varepsilon_{\vk S}(T_k)\mathfrak{C}_k(T_k)}
}
takes place, where
\bqny{
\psi^*(\vk S)&=&\pk{\exists_{t\in[0,T_k]}:\vk Z^{*k}(t,\vk t^k)\in\vk S~for~some~\vk t^k\in\mathbb{T}_{k-1}},\\
\nu(\vk S)&=&\pk{\vk Z^{*k}(T_k,\vk t^k)\in\vk S~for~some~\vk t^k\in\mathbb{T}_{k-1}}.
}
Define the following stopping time:
\bqny{
\tau_k=\inf\left\{t: \vk Z^{*k}(t,\vk t^k)\in u\vk S~for~some~\vk t^k\in\mathbb{T}^k\right\},
}
and the random vector
\bqny{
\tilde{\vk x}_k=\begin{cases}
\vk x^*, \qquad &\tau_k\leq T_k,\\
\vk 0,\qquad &otherwise,
\end{cases}
}
where $\vk x^*$ is any point from the following set: 
$$
\bigcup\limits_{\vk t^k\in\mathbb{T}^k}\left\{\vk Z^{*k}(\tau_k,\vk t^k)\right\}\bigcap u\vk S.
$$. 
Using the total probability formula we obtain
\bqny{
\nu(u\vk S)=\int_0^{T_k}\int_{u\partial \vk S}\pk{\tilde{\vk x}_k\in\td\vk x_0,\tau_k\in\td t}\pk{\vk Z^{*k}(T_k,\vk t^k)\in u\vk S~for~some~\vk t^k\in\mathbb{T}_{k-1}\left|\tau_{k}=t,\tilde{\vk x}_k=\vk x_0\right.}.
}
For any $\vk t^k\in\mathbb{T}_{k-1}$ we have
\bqny{
\vk Z^{*k}(T_k,\vk t^k)-\vk Z^{*k}(t,\vk t^k)=\vk Z_k(T_k)-\vk Z_k(t)-(\vk c_k(T_k)-\vk c_k(t)).
}
Thus, using the same chain of inequalities as in \netheo{ConeSet} we obtain
\bqny{
& &\pk{\vk Z^{*k}(T_k,\vk t^k)\in u\vk S~for~some~\vk t^k\in\mathbb{T}^k\left|\tau_{k}=t,\tilde{\vk x}_k=\vk x_0\right.}\\
& &\qquad\geq\pk{\vk Z_k(T_k)-\vk Z_k(t)-(\vk c_k(T_k)-\vk c_k(t))\in u\vk S-\vk x_0}\\
& &\qquad \geq \frac{\mathfrak{C}_k(T_k)\varepsilon_{\vk S}(T_k)}{2^{d/2}},
}
which completes the proof.
\QED

\begin{remark}
	The random variable $\tau_k$ is measurable, because it can be represented as
	\bqny{
		\tau_{k}=\inf\left\{t:~\vk Z_k(t)-\vk c_k(t)\in \vk S^*_k \right\},
	}
	where
	\bqny{
		\vk S^*_k=\bigcup\limits_{t^k\in \mathbb{T}^k} \left(u\vk S-\sum_{i=1}^{k-1}\left(\vk x_i(t_i)-\vk c_i(t_i)\right)-\sum_{i=k+1}^{n}\left(\vk x_i-\vk c_i(T)\right)\right).
	}
	As all the functions $\vk x_i(t)$ and $\vk c_i(t)$ are continuous, set $\vk S_k^*$ is closed, hence $\tau_k$ is measurable.
\end{remark}

\prooftheo{TimeTransform}
Define a stopping the
\bqny{
\tau=\inf\{t\in [0,T]:\vk Z(\vk f(t))-\vk c(t) \in u\vk S \}.
}
According to the strong Markov property
\bqny{
\pk{\vk Z(\vk f(T)) -\vk c (T)\in u\vk S}&=&\int_{0}^{T}\int_{u\partial \vk S}\pk{\vk Z(\vk f(\tau))-\vk c(\tau)\in\td \vk x, \tau\in \td t}\\
& &\qquad\qquad\times\pk{\vk Z(\vk f(T))-\vk c(T)\in u\vk S|\vk Z(\vk f(t))-\vk c(t)=\vk x,\tau=t}.
}
In view of \nelem{1}, $u\vk S$ satisfies the cone condition with respect to the process $\vk Z(t)$. Consequently, we have
\bqny{
\lefteqn{\pk{\vk Z(\vk f(T))-\vk c(T)\in u\vk S|\vk Z(\vk f(t))-\vk c(t)=\vk x,\tau=t}}\\
& &=\pk{\vk Z(\vk f(T))-\vk c(T)\in u\vk S|\vk Z(\vk f(t))-\vk c(t)=\vk x}\\
& &\geq\pk{\vk Z(\vk f(T))-\vk c(T)\in V_{\vk x}|\vk Z(\vk f(t))-\vk c(t)=\vk x}\\
& &=\pk{\vk Z(\vk f(T))-\vk Z(\vk f(t))-\vk c(T)+\vk c(t)\in V_{\vk x}-\vk x}\\
& &=\pk{\vk Z(\vk f(T)-\vk f(t))-(\vk c(T)-\vk c(t))\in V_{\vk x}-\vk x}\\
& &=\pk{\vk Z(\vk\delta(t))-\frac{\vk c(T)-\vk c(t)}{\sqrt{f_1(T)-f_1(t)}}\in \frac{V_{\vk x}-\vk x}{\sqrt{f_1(T)-f_1(t)}}}\\
& &\geq\pk{\vk Z(\vk\delta(t))-\frac{\vk c(T)-\vk c(t)}{\sqrt{f_1(T)-f_1(t)}}\in \frac{V_{\vk x}-\vk x}{\sqrt{f_1(T)}}}\\
& &=\int_{\vk y\in\frac{V_{\vk x}-\vk x}{\sqrt{f_1(T)}}}\varphi_{\vk\delta(t)}\left(\vk y+\frac{\vk c(T)-\vk c(t)}{\sqrt{f_1(T)-f_1(t)}}\right)\td \vk y\\
& &\geq\int_{\vk y\in\frac{V_{\vk x}-\vk x}{\sqrt{f_1(T)}}}\mathfrak{C}(T)\varphi_{\vk\delta(t)}(\sqrt{2}\vk y)\td \vk y\\
& &\geq\frac{\mathfrak{C}(T)}{2^{d/2}}\pk{\vk Z(\vk \delta(t))\in\frac{V_{\vk x}-\vk x}{\sqrt{f_1(T)}}}\\
& &\geq\frac{\mathfrak{C}(T)}{(2f_1(T))^{d/2}}\pk{\vk Z(\vk \delta(t))\in V_{\vk x}-\vk x}\\
& &=\frac{\mathfrak{C}(T)}{(2f_1(T))^{d/2}}\pk{\vk B(\vk \delta(t))\in A^{-1}(V_{\vk x}-\vk x)}\\
& &=\frac{\mathfrak{C}(T)}{(2f_1(T))^{d/2}}\frac{1}{\sqrt{2\pi\prod_{i=1}^d\delta_i(t)}}\int_{\vk y\in A^{-1}(V_{\vk x}-\vk x)} e^{-\frac{1}{2}\sum\limits_{i=1}^{d}\frac{y_i^2}{\delta_i(t)}}\td\vk y,
}
where $\varphi_{\vk\delta(t)}$ is the pdf of $\vk Z(\vk \delta(t))$. Using that all the functions $\delta_i(t)$ are bounded and separated from zero for $t\in [0,T]$, there exists some constants $\delta,\Delta>0$, such that for all $i\in\{1\ldot d\}$ and all $t\in [0,T]$ 
$$
\delta\leq\delta_i(t)\leq \Delta.
$$
Hence we obtain
\bqny{
\frac{1}{\sqrt{2\pi\prod_{i=1}^d\delta_i(t)}}\geq \frac{1}{\sqrt{2\pi\prod_{i=1}^d\Delta}},\qquad e^{-\frac{1}{2}\sum\limits_{i=1}^{d}\frac{x_i^2}{\delta_i(t)}}\geq e^{-\frac{1}{2}\sum\limits_{i=1}^{d}\frac{x_i^2}{\delta}},
}
and finally
\bqny{
\lefteqn{\pk{\vk Z(\vk f(T))-\vk c(T)\in u\vk S|\vk Z(\vk f(t))-\vk c(t)=x,\tau=t}}\\
& &\geq \frac{\mathfrak{C}(T)}{(2f_1(T))^{d/2}}\frac{1}{\sqrt{2\pi\prod_{i=1}^d\Delta}}\int_{\vk y\in A^{-1}(V_{\vk x}-\vk x)} e^{-\frac{1}{2}\sum\limits_{i=1}^{d}\frac{y_i^2}{\delta}}\td\vk y\\
& &= \frac{\mathfrak{C}(T)}{(2f_1(T))^{d/2}}\frac{\sqrt{\prod_{i=1}^d\delta}}{\sqrt{\prod_{i=1}^d\Delta}}\pk{\vk B(\delta)\in A^{-1}(V_{\vk x}-\vk x)}\\
& &\geq \frac{\mathfrak{C}(T)}{(2f_1(T))^{d/2}}\frac{\sqrt{\prod_{i=1}^d\delta}}{\sqrt{\prod_{i=1}^d\Delta}}\varepsilon_{\vk S}(\delta).
}
Hence the claim follows.
\QED 

\vspace{1cm}

\begin{center}
	{\bf Acknowledgements}
\end{center}  

I am grateful to four reviewers for numerous comments and  suggestions  that lead to a significant improvement of the manuscript.
Partial financial support from the SNSF Grant 200021-191984 is kindly acknowledged.

\bibliographystyle{ieeetr}
\bibliography{EEEA}

\end{document}